\documentstyle[12pt]{article}
\if@twoside  m
\else
\fi
\newcommand{\ie}{{\em i.e.}}
\newcommand{\HH}{{\cal H}}
\newcommand{\OO}{{\cal O}}
\newcommand{\R}{I\!\!R}
\newcommand{\im}{{\rm Im\,}}
\newcommand{\QED}{\mbox{\rule[-1.5pt]{6pt}{10pt}}}
\newtheorem{claim}{Claim}
\newtheorem{proposition}[claim]{Proposition}
\newtheorem{theorem}[claim]{Theorem}
\newtheorem{lemma}[claim]{Lemma} \newtheorem{remark}[claim]{Remark}
 
\begin{document}
\vspace*{20mm}  \noindent
{\Large\bf Weakly coupled states on branching graphs}
\vspace{5mm}  \noindent
{\bf Pavel Exner} \\ {\small
Nuclear Physics Institute, Academy of Sciences, 25068 \v{R}e\v{z}
near Prague, \\ and Doppler Institute, Czech Technical University,
B\v rehov\'a 7, 11519 Prague, \\ Czech Republic, \\ \em exner@ujf.cas.cz}
\vspace{10mm}  \noindent {\small {\bf Abstract.} We consider a Schr\"odinger
particle on a
graph consisting of $\,N\,$ links joined at a single point. Each link
supports a real locally integrable potential $\,V_j\,$; the
self--adjointness is ensured by the $\,\delta\,$ type boundary
condition at the vertex. If all the links are semiinfinite and
ideally coupled, the potential decays as $\,x^{-1-\epsilon}$ along
each of them, is non--repulsive in the mean and weak enough, the
corresponding Schr\"odinger operator has a single negative
eigenvalue; we find its asymptotic behavior. We also derive a bound
on the number of bound states and explain how the $\,\delta\,$
coupling constant may be interpreted in terms of a family of squeezed
potentials.}  \vspace{10mm}  \noindent Recent progress in investigation of
``mesoscopic" systems attracted a
wave of attention to properties of quantum mechanical particles whose
motion is confined to a graph --- see \cite{Ad,AL,ARZ,BuT,ES,GPS,GLR,GP} and
references therein. The problem is not new; it appeared for the first time in
early fifties
in connection with the free--electron model of organic molecules
\cite{RuS}. However, the mentioned studies brought not only physical
applications but fresh mathematical insights as well.  They concern, in
particular, relations between spectral properties
and the dimensionality of the configuration space. For instance, band
structure of the spectrum for a periodic rectangular--lattice graph
exhibits an interesting dependence on number--theoretic properties of
the ratio $\,\theta\,$ of the rectangle sides \cite{E2,E3,EG}. For
most values of $\,\theta\,$ there are infinitely many gaps once the
the ``coupling constant" at graph vertices is nonzero; hence such
systems do not conform with the Bethe--Sommerfeld conjecture \cite{Sk1,Sk2} and
behave rather as one--dimensional ones. On the other hand, for a
``bad--irrational" lattice there are no gaps if the
coupling at the vertices is weak enough which is {\em not} a typical
one--dimensional behavior; infinitely many gaps open only above a
certain critical value.  In the present letter we are going to show that the
mentioned
one--dimensional feature is contained already in the way in which the
wavefunctions are coupled at graph vertices as long as the number of
links entering a vertex is finite. In fact, we shall demonstrate that
{\em continuous} Schr\"odinger operators on graphs may be regarded as
an extension of the standard Sturm--Liouville theory of second--order
ODE's.  For simplicity, we restrict ourselves to the simplest situation of a
graph with a single vertex in which a finite number of $\,N\ge 2\,$
links of lengths $\,\ell_j\,,\; j=1,\dots,N\,$ are joined. Since each
link can be mapped to a finite or semiinfinite interval, the
corresponding state Hilbert space is identified with
$\,\HH:=\bigoplus_{j=1}^N L^2(0,\ell_j)\,$. We suppose that the
$\,j$th--link supports a real--valued potential $\,V_j\,$ and assume
that
\begin{description} \item{\em (p1)} $\;V_j \in L^1_{loc}(0,\ell_j)\,,\;\;
j=1,\dots,N\,$, \item{\em (p2)} if a given link is semiinfinite, the
differential expression $\,-{d^2\over dx^2}+ V_j(x)\,$ is LPC at infinity; if
$\,\ell_j<\infty\,$ this requirement is replaced by a fixed boundary
condition, $\,\psi(\ell_j)\cos\omega_j+ \psi'(\ell_j)\sin\omega_j
=0\,$.
\end{description}
Given a family $\,V:=\{V_j\}_{j=1}^N\,$ and $\,\alpha\in\R \cup\{\infty\}\,$,
we define the operator $\,H_{\alpha}(V)\,$ by     \begin{equation}
\label{Hamiltonian}
H_{\alpha}(V)\{\psi_j\}\,:=\,
\{-\psi_j''+V_j\psi_j\}
\end{equation}
with the natural domain requirements at each link and the boundary
conditions
\begin{equation} \label{bc} \psi_1=\cdots=\psi_n=:\psi\,, \qquad \sum_{j=1}^N
\psi'_j= \alpha\psi     \end{equation}
at the vertex, where $\,\psi_j:= \lim_{x\to 0+}
\psi_j(x)\,$ and $\,\psi'_j:=\lim_{x\to 0+} \psi'_j(x)\,$  are the
corresponding boundary values. For $\,\alpha=\infty\,$ the
requirement (\ref{bc}) is replaced by the Dirichlet condition,
$\,\psi_j=0\,,\; j=1,\dots,N\,$; in that case the operator is
decoupled, $\,H_{\infty}(V)= \bigoplus_{j=1}^N h_j(V_j)\,$.  The
condition (\ref{bc}) is known to produce a self--adjoint operator if
$\,V=0\;$ \cite{ES}; it is straightforward to check that under the
stated assumptions this conclusion is not changed:

\begin{proposition}
The operator $\,H_{\alpha}(V)\,$ is self--adjoint for any
$\,\alpha\in\R \cup\{\infty\}\,$ and $\,V\,$ obeying the conditions
(p1,2).
\end{proposition}  Our first result consists of showing that if the coupling is
ideal,
$\,\alpha=0\,$, and all the links are semiinfinite, an arbitrarily
weak potential potential which is not repulsive in the mean and
decays fast enough produces a bound state.  \begin{theorem}
\label{weak coupling thm} Let $\,\ell_j=\infty\,$ and $\,V_j\in
L^2(\R_+,(1+|x|)dx)\,,\; j=1,\dots,N\,$. Then the operator
$\,H_0(\lambda V)\,$ has for all sufficiently small $\,\lambda>0\,$ a
single negative eigenvalue $\,\epsilon(\lambda)=
-\kappa(\lambda)^2\;$ {\em iff}
\begin{equation} \label{potential mean} \sum_{j=1}^N \int_0^{\infty} V_j(x)\,dx
\,\le\,0\,.    \end{equation}
In that case, its asymptotic behavior is given by \begin{eqnarray} \label{weak
coupling} \kappa(\lambda) \! &=& \! -\,{\lambda\over N}\,
\sum_{j=1}^N \int_0^{\infty} V_j(x)\,dx \,-\, {\lambda^2\over
2N}\, \Biggl\lbrace\, \sum_{j=1}^N \int_0^{\infty}\! \int_0^{\infty}
V_j(x) |x\!-\!y| V_j(y)\,dx\,dy \nonumber \\ \\ &+& \! \sum_{j,\ell=1}^N
\left( {2\over N}\,-\,\delta_{j\ell}\right)\, \int_0^{\infty}\!
\int_0^{\infty} V_j(x) (x\!+\!y) V_{\ell}(y)\, dx\,dy\, \Biggl\rbrace
\,+\, \OO(\lambda^3)\,. \nonumber \end{eqnarray}    \end{theorem}  Under the
assumptions, the essential spectrum of the Dirichlet link
operators $\,h_j(V_j)\,$ is $\,[0,\infty)\,$. As we shall show below,
$\,H_{\alpha}(V)\,$ and $\,H_{\infty}(V)\,$ differ by a rank--one
perturbation in the resolvent, hence we have also
$\,\sigma_{ess}(H_{\alpha}(V))= [0,\infty)\,$. Moreover, the
operators $\,h_j(\lambda V_j)\,$ have no discrete spectrum for small
$\,\lambda\,$, so $\,H_{\alpha}(\lambda V)\,$ has at most one
negative eigenvalue.  To prove the existence condition (\ref{potential mean})
and the
asymptotic expansion (\ref{weak coupling}), we employ the explicit
form of the resolvent; we shall derive it for the general case.
Denote $\,H_{\alpha}:= H_{\alpha}(V)\,$. Given $\,k\,$ with $\,\im
k\ge 0\,, \;k^2\in \rho(H_{\infty})\,$, we denote by $\,u_j\equiv
u_j(\cdot;k)\,$  and $\,v_j\equiv v_j(\cdot;k)\,$ solutions to
$\,-\psi''_j+V_j\psi_j=
k^2\psi_j\,$ with the appropriate behavior at $\,0\,$ and
$\,\ell_j\,$, respectively, \ie, $\,u_j(0;k)=0\,$ while $\,v_j\,$ is
square integrable at infinity if $\,\ell_j=\infty\,$ and satisfies
the fixed boundary condition otherwise. If $\,\alpha=\infty\,$, the
links are decoupled and the corresponding components of
$\,H_{\infty}\,$ are characterized by the resolvent kernels
\begin{equation} \label{link resolvent kernel}
g_j(x,y;k)\,:=\,
-\,{u_j(x_<;k)v_j(x_>;k)\over W(u_j,v_j)}\,,
\end{equation}
where conventionally $\,x_<:=\min\{x,y\}\,$ and similarly for
$\,x_>\,$, and $\,W(u_j,v_j)\,$ is the Wronskian of the two
solutions.

\begin{lemma} \label{resolvent kernel lemma}
Assume (p1,2). For arbitrary $\,\alpha\in\R\,$ and $\,k\,$ such that
$\,\im k\ge 0\,$ and $\,k^2\in \rho(H_{\infty})\cap \rho(H_{\alpha})\,$, the
resolvent $\,(H_{\alpha}\!-k^2)^{-1}$ is a matrix integral operator with the
kernel \begin{equation}
\label{resolvent kernel} G^{\alpha}_{j\ell}(x,y;k)\,=\,
\delta_{j\ell} g_j(x,y;k)\,+\, {v_j(x;k)v_{\ell}(y;k)\over
v_j(0;k)v_{\ell}(0;k) (\alpha-M(k))}\,, \end{equation} where
$\,M(k):= \sum_{j=1}^N {v_j'(0;k)\over v_j(0;k)}\,$.
\end{lemma}
{\em Proof:} By Krein's formula \cite[Appendix A]{AGHH}, the sought
kernel is of the form $\,\delta_{j\ell} g_j(x,y;k)+ \Lambda_{j\ell}
v_j(x;k)v_{\ell}(y;k)\,$. To find the unknown coefficients, we
express from here $\,\psi:= (H_{\alpha}\!-k^2)^{-1}\phi\,$. Since
$\,\psi\,$ has to satisfy the boundary conditions (\ref{bc}) for any
$\,\phi\in\HH\,$, we arrive at a system of linear equations for
$\,\Lambda_{j\ell}\,$ which yields (\ref{resolvent kernel}). \quad \QED 
\begin{remark}
{\rm Eigenvalues of $\,H_{\alpha}(V)\,$ are determined by zeros of
the denominator in (\ref{resolvent kernel}).  Since all the
logarithmic derivatives are decreasing functions of $\,k\,$, so is
$\,M(k)\,$, and each interval between neigbouring points of
$\,\bigcup_{j=1}^N \sigma_{disc}(h_j(V_j))\,$ contains for an
arbitrary $\,\alpha \in\R\,$ just one simple eigenvalue. Hence the
multiplicity of $\,\sigma_{disc}(H_{\alpha}(V))\,$ in the coupled
case, $\,\alpha\ne\infty\,$, does not exceed $\,N\!-\!1\,$; in
particular, the discrete spectrum is simple for $\,N=2\,$. On the
other hand, eigenvalues of multiplicity $\,N\!-\!1\,$ arise naturally
if $\,H_{\infty}(V)\,$ has a ``fully degenerate" eigenvalue as
noticed in \cite{GPS}. }
\end{remark}
{\em Proof of Theorem~\ref{weak coupling thm}:} We have
$\,H_0(\lambda V)= A_0\!+\lambda V\,$, where $\,A_0:=H_0(0)\,$. The
well--known resolvent formula \cite[Appendix B]{AGHH} shows that
possible negative eigenvalues of $\,H_0(\lambda V)\,$ are given by
the Birman--Schwinger principle: such an eigenvalue exists {\em iff}
$\;\lambda K\,$ has the eigenvalue $\,-1\,$, where $\,K:=
|V|^{1/2}(A_0\!-k^2)^{-1} V^{1/2}\,$ is determined by its kernel,
\begin{eqnarray} K_{j\ell}(x,y;k) \! &:=&\!
\delta_{j\ell} |V_j(x)|^{1/2} {\sinh \kappa x_<\, e^{-\kappa x_>}\over
\kappa}\, V_j(y)^{1/2} \nonumber \\ \\ &+& \! {1\over \kappa N}\,
|V_j(x)|^{1/2} e^{-\kappa(x+y)} V_{\ell}(y)^{1/2}\,, \nonumber
\end{eqnarray}
where we have introduced $\,\kappa:=-ik\,$ and $\,V_j(\cdot)^{1/2}$
is the standard signed square root of the potential. Writing $\,K=
P_{\kappa}\!+Q_{\kappa}\,$, we find that the first part can be
estimated as
$$ \|P_{\kappa}\|^2\,\le\, \|P_{\kappa}\|^2_{HS}\,\le\,  \sum_{j=1}^N
\int_0^{\infty} \int_0^{\infty} |V_j(x)| x_<^2 |V_j(y)|\, dx\,dy \,\le\,
\sum_{j=1}^N \left(\int_0^{\infty}
x|V_j(x)|\, dx \right)^2\!,
$$
so it has a bound independent of $\,\kappa\,$ and the operator
$\,P_0\,$ corresponding to the kernel $\,P_{0,j\ell}(x,y;k):= \delta_{j\ell}
V_j(x)\, x_<^2 V_j(y)\,$ is well--defined. Furthermore, $\,Q_{\kappa}\,$ is
rank--one and the same is true for
\begin{equation} \label{1D part} \lambda(I+\lambda
P_{\kappa})^{-1}Q_{\kappa}\,=\, (\psi,\cdot)\phi\,,    \end{equation}
where
$$ \psi\,:=\, {\lambda\over \kappa N}\,e^{-\kappa\cdot} V^{1/2}\,, \quad
\phi\,:=\, (I+\lambda P_{\kappa})^{-1} e^{-\kappa\cdot} |V|^{1/2} \;;
$$
due to the above result the inverse makes sense for
all $\,\lambda\,$ small enough. Using the identity
$$
(I+\lambda K)^{-1}\,=\, [I+\lambda(I+\lambda P_{\kappa})^{-1}
Q_{\kappa}]^{-1} (I+\lambda P_{\kappa})^{-1}
$$
we find that for small $\,\lambda\,$ eigenvalues of $\,\lambda K\,$
coincide with those of the operator (\ref{1D part}). However, the
latter has just one eigenvalue $\,\xi(\lambda)= (\psi,\phi)\,$.
Demanding $\,\xi(\lambda)= -1\,$ we get an implicit equation for
$\,\kappa(\lambda)\,$; the rest of the argument is the same as
in \cite{BGS,Kl}. If the potential mean (\ref{potential mean}) equals  zero,
the linear term is absent as well as the $\,2/N\,$ part of the
quadratic one, so
$$ \kappa(\lambda)\,=\,-\, {\lambda^2\over 4N}\, \sum_{j=1}^N \int_{\R^2}\!
\tilde V_j(x)\,|x\!-\!y|  \tilde V_j(y)\,dx\,dy \,+\,\OO(\lambda^3)\,, $$
where $\,\tilde V_j\,$ is the odd extension of $\,V_j\,$ to the
whole $\,\R\,$; hence the standard argument works again. \quad
\QED \vspace{3mm}  In the case $\,N=2\,$ the asymtotic expansion (\ref{weak
coupling})
reduces to the well--known formula for Schr\"odinger operators on
line \cite{BGS,Kl}; the present proof illustrates why the
factorization of the singular part in these papers (in contrast to
\cite{Si}) was optimal. In the same way one can generalize the
trace--class bound on the number of one--dimensional Schr\"odinger
operators \cite{Se,Kl,Ne}.

\begin{proposition} \label{bound-state number proposition}
Let $\,\ell_j=\infty\,$ and
$\,V_j^{(-)}:=\max\{0,-V_j\}\in L^2(\R_+,(1+|x|)dx)\,$ for
$\,j=1,\dots,N\,$. Then the number of negative eigenvalues of
$\,H_0(V)\,$ may be estimated by
\begin{eqnarray} \label{bound-state number}
N(V^{(-)})\!& \le &\! 1\,+\, {1\over 2\langle
V^{(-)}\rangle}\, \Biggl\{\sum_{j=1}^N \int_0^{\infty}\! \int_0^{\infty}
|x\!-\!y|\, V_j^{(-)}(x)V_j^{(-)}(y)\,dx\,dy \nonumber \\ \\ \!& +
&\!\sum_{j,\ell=1}^N \left( {2\over N}\,-\, \delta_{j\ell}\right)\,
\int_0^{\infty}\! \int_0^{\infty}  (x\!+\!y)\,V_j^{(-)}(x)V_{\ell}^{(-)}(y)\,
dx\,dy \,\Biggr\}\,, \nonumber    \end{eqnarray}
where $\,\langle V^{(-)}\rangle:= \sum_{j=1}^N
\int_0^{\infty} V_j^{(-)}(x)\,dx\,$.
\end{proposition}  It is well known that the $\,\delta\,$ interaction on line
can be
approximated by means of a family of squeezed potentials. In the
present context we can extend this result to branching graphs; this
gives the boundary condition (\ref{bc}) with a nonzero $\,\alpha\,$ a
natural meaning of low--energy description of a non--ideal junction.
Given $\,W_j\,$ we define the scaled potentials by
\begin{equation}
\label{scaled potentials} W_{\epsilon,j}\,:=\, {1\over\epsilon}\,
W_j\left(x\over\epsilon \right)\,,\quad j=1,\dots,N\,.
\end{equation}     \begin{theorem} \label{approximation thm}
Suppose that $\,V_j \in L^1_{loc}(0,\ell_j)\,$ are below bounded and
$\,W_j\in L^1(0,\ell_j)\,$ for $\,j=1,\dots,N\,$. Then
\begin{equation}
\label{approximation} H_0(V+W_{\epsilon})\,\longrightarrow\,
H_{\alpha}(V) \qquad {\rm as} \quad \epsilon\to 0+      \end{equation}
in the norm resolvent sense, where $\,\alpha= \langle W\rangle:= \sum_{j=1}^N
\int_0^{\infty} W_j(x)\,dx\,$.    \end{theorem}
{\em Proof:} Let $\,G_{j\ell}^{W_\epsilon}(x,y;k)\,$
denote the resolvent kernel of $\,H_0(V+W_{\epsilon})\,$. Using once
more the resolvent formula, we may rewrite it as
\begin{eqnarray*}
G_{j\ell}^{W_\epsilon}(x,y;k) \!&=&\! G_{j\ell}^0(x,y;k) \,-\, \sum_{r,s}\,
\int_0^{\infty} \int_0^{\infty} G_{j\ell}^0(x,y';k)\,  W_{\epsilon,r}(x')^{1/2}
 \\ \\ & \times &\! \left(
I+|W_{\epsilon}|^{1/2} \left(H_0(V)-k^2\right)^{-1}
W_{\epsilon}^{1/2}\right)^{-1}_{rs}(x',x'')\,
|W_{\epsilon,r}(x'')|^{1/2} \\ \\ & \times &\!
G_{s\ell}^0(x'',y;k)\,dx'\,dx''\,.  \end{eqnarray*} Changing the integration
variables to $\,x'/\epsilon\,$ and
$\,x''/\epsilon\,$ as in \cite[Sec.I.3.2]{AGHH}, we can rewrite the
resolvent in question in the form $\,-B_{k,\epsilon}\, \left(
I+C_{k,\epsilon}\right)^{-1} \tilde B_{k,\epsilon}\,$, where the involved
operators are determined by their kernels,
\begin{eqnarray*}
(B_{k,\epsilon})_{j\ell}(x,y) \!&=&\!
G_{j\ell}^0(x,\epsilon y;k) W_{\ell}(y)^{1/2}\,, \\ \\ (\tilde
B_{k,\epsilon})_{j\ell}(x,y) \!&=&\! |W_j(x)|^{1/2}
G_{j\ell}^0(\epsilon x,y;k) \,, \\ \\ (C_{k,\epsilon})_{j\ell}(x,y) \!&=&\!
|W_j(x)|^{1/2}
G_{j\ell}^0(\epsilon x,\epsilon y;k)\, W_{\ell}(y)^{1/2}\,,
\end{eqnarray*} which converge pointwise to \begin{eqnarray*}
(B_k)_{j\ell}(x,y) \!&=&\! G_{j\ell}^0(x,0;k) W_{\ell}(y)^{1/2}\,, \\
\\ (\tilde B_k)_{j\ell}(x,y) \!&=&\! |W_j(x)|^{1/2}
G_{j\ell}^0(0,y;k) \,, \\ \\ (C_k)_{j\ell}(x,y) \!&=&\! |W_j(x)|^{1/2}
G_{j\ell}^0(0,0;k)\,
W_{\ell}(y)^{1/2}\,,
\end{eqnarray*}
respectively, as $\,\epsilon\to 0+\,$. The explicit form of the last
operator makes it possible to find the inverse,
$$
(I+C_k)^{-1}_{j\ell}(x,y)\,=\, \delta(x\!-\!y)
\delta_{j\ell} \,-\, {|W_j(x)|^{1/2} W_{\ell}(y)^{1/2}\over \langle W\rangle
 -M(k)}\,,
$$
where we have employed (\ref{resolvent kernel}). The
correction to the ``free" resolvent kernel is therefore equal to
\begin{eqnarray*}
&& -\,\sum_{r}\, \int_0^{\infty}\,dx'\, W_r(x')\,
G_{jr}^0(x,0;k) G_{r\ell}^0(0,y;k) \\ \\ && +\, \sum_{r,s}\,
\int_0^{\infty}\, \int_0^{\infty}\, dx'\, dx''\, W_r(x')W_s(x'')\,
 {G_{jr}^0(x,0;k) G_{r\ell}^0(0,y;k)\over \langle W\rangle -M(k)}  \\
\\ && =\, {v_j(x;k) v_{\ell}(y;k)\over v_j(0;k) v_{\ell}(0;k)}\:
{\langle W\rangle\over M(k)(\langle W\rangle-M(k))}\;;
\end{eqnarray*}
this coincides with the analogous term in (\ref{resolvent kernel}) if
we set $\,\alpha:= \langle W\rangle\,$.  It is sufficient to check
the norm--resolvent convergence for a particular $\,k\,$.  Since all
the $\,h_j(V_j)\,$ are below bounded by assumption, one may choose
$\,\kappa\,$ large enough to get exponentially decaying solutions
$\,v_j(\cdot;k)\,$ at the semiinfinite links; the argument then
proceeds as in \cite[Sec.I.3.2]{AGHH}. \quad \QED \vspace{10mm}  \noindent {\em
Acknowledgments.} The author is grateful for the hospitality
extended to him in the Centre de Physique Th\'eorique, C.N.R.S.,
Marseille--Luminy, where this work was done. The research has been
partially supported by the Grant AS
No.148409.    
\vspace{5mm} \begin{thebibliography}{article}    \bibitem[Ad]{Ad}
V.M.~Adamyan: Scattering matrices for
microschemes, {\em Oper.Theory: Adv. Appl.} {\bf 59} (1992), 1--10.
\vspace{-1.8ex}
\bibitem[AGHH]{AGHH} S.~Albeverio, F.~Gesztesy,
R.~H\o egh-Krohn, H.~Holden: {\em Solvable Models in Quantum
Mechanics}, Springer, Heidelberg 1988.  \vspace{-1.8ex}
\bibitem[AL]{AL} Y.~Avishai, J.M.~Luck: Quantum percolation and
ballistic conductance on a lattice of wires, {\em Phys.Rev.} {\bf
B45} (1992), 1074--1095.  \vspace{-1.8ex}
\bibitem[ARZ]{ARZ}
J.E.~Avron, A.~Raveh, B.~Zur: Adiabatic transport in multiply
connected systems, {\em Rev.Mod.Phys.} {\bf 60} (1988), 873--915.
\vspace{-1.8ex}
\bibitem[BGS]{BGS} R.~Blanckenbecler,
M.L.~Goldberger, B.~Simon: The bound states of weakly coupled
long-range one-dimensional quantum Hamiltonians, {\em Ann.Phys.} {\bf
108} (1977), 89-78.  \vspace{-1.8ex}
\bibitem[BT]{BuT} W.~Bulla,
T.~Trenckler: The free Dirac operator on compact and non-compact
graphs, {\em J.Math.Phys.} {\bf 31} (1990), 1157-1163.
\vspace{-1.8ex}
\bibitem[E2]{E2} P.~Exner: Lattice Kronig--Penney
models, {\em Phys.Rev.Lett.} {\bf 74} (1995), 3503--3506.
\vspace{-1.8ex}
\bibitem[E3]{E3} P.~Exner: Contact interactions on
graph superlattices, {\em J.Phys.} {\bf A}, to appear \vspace{-4.8ex}
\bibitem[EG]{EG} P.~Exner, R.~Gawlista: Band spectra of rectangular
graph superlattices, {\em Phys.Rev.} {\bf B}, to appear
\vspace{-1.8ex}
\bibitem[E\v S]{ES} P.~Exner, P.~\v{S}eba: Free
quantum motion on a branching graph, {\em Rep. Math.Phys.} {\bf 28}
(1989), 7--26.  \vspace{-1.8ex}
\bibitem[GPS]{GPS} A.~Gangopadhyaya,
A.Pagnamenta, U.~Sukhatme: Quantum mechanics of multi--prong
potentials, {\em J.Phys.} {\bf A28} (1995), 5331--5347.
\vspace{-1.8ex}
\bibitem[GP]{GP} N.I.~Gerasimenko, B.S.~Pavlov:
Scattering problem on noncompact graphs, {\em Teor.Mat.Fiz.} {\bf 74}
(1988), 345--359.  \vspace{-1.8ex}
\bibitem[GLRT]{GLR} J.~Gratus,
C.J.~Lambert, S.J.~Robinson, R.W.~Tucker: Quantum mechanics on
graphs, {\em J.Phys.} {\bf A27} (1994), 6881--6892.  \vspace{-1.8ex}
\bibitem[Kl]{Kl} M.~Klaus: On the bound state of Schr\"odinger
operators in one dimension, {\em Ann.Phys.} {\bf 108} (1977),
288--300.  \vspace{-1.8ex}
\bibitem[Ne]{Ne} R.G.~Newton: Bounds on
the number of bound states for the Schr\"odinger equation in one and
two dimension, {\em J. Operator Theory} {\bf 10} (1983), 119--125.
\vspace{-4.8ex}
\bibitem[RuS]{RuS} K.~Ruedenberg, C.W.~Scherr:
Free--electron network model for conjugated systems, I.~Theory, {\em
J.Chem.Phys.} {\bf 21} (1953), 1565--1581.  \vspace{-1.8ex}
\bibitem[Se]{Se} N.~Seto: Bargmann's inequality in spaces of
arbitrary dimension, {\em Publ. RIMS} {\bf 9} (1974), 429--461.
\vspace{-1.8ex}
\bibitem[Si]{Si} B.~Simon: The bound state of weakly
coupled Schr\"odinger operators in one and two dimensions, {\em
Ann.Phys.} {\bf 97} (1976), 279-288.  \vspace{-1.8ex}
\bibitem[Sk1]{Sk1} M.M.~Skriganov: Proof of the Bethe--Sommerfeld
conjecture in dimension two, {\em Sov.Math.Doklady} {\bf 20} (1979),
956-959.  \vspace{-1.8ex}
\bibitem[Sk2]{Sk2} M.M.~Skriganov: The
multidimensional Schr\"odinger operator with a periodic potential,
{\em Math.USSR Izvest.} {\bf 22} (1984), 619--645.
\vspace{-1.8ex} \end{thebibliography}
\end{document}